\documentclass[a4paper,11pt,reqno]{amsart}

\usepackage{amssymb}

\theoremstyle{definition}

\newtheorem*{example*}{Example}

\newtheorem*{remark*}{Remark}

\newcommand{\frg}{{\mathfrak g}}

\newcommand{\calL}{{\mathcal L}}

\providecommand{\espan}[1]{\text{span}\left\{ #1\right\}}

 \DeclareMathOperator{\End}{End}

\input prepictex
\input pictex
\input postpictex

\begin{document}

\title{A Lie grading which is not a semigroup grading}

\author[Alberto Elduque]{Alberto Elduque$^{\star}$}
 \thanks{$^{\star}$ Supported by the Spanish Ministerio de
 Educaci\'{o}n y Ciencia
 and FEDER (MTM 2004-081159-C04-02) and by the
Diputaci\'on General de Arag\'on (Grupo de Investigaci\'on de
\'Algebra)}
 \address{Departamento de Matem\'aticas, Universidad de
Zaragoza, 50009 Zaragoza, Spain}
 \email{elduque@unizar.es}

\date{December 15, 2005}



%

\maketitle



Patera and Zassenhaus \cite{PZ} define a \emph{Lie grading} as a
decomposition of a Lie algebra into a direct sum of subspaces
\begin{equation}\label{eq:Liegr}
\calL=\oplus_{g\in G}\calL_g,
\end{equation}
such that $\calL_g\ne 0$ for any $g\in G$, and for any $g,g'\in G$,
either $[\calL_g,\calL_{g'}]=0$ or there exists a $g''\in G$ such
that $0\ne [\calL_g,\calL_{g'}]\subseteq \calL_{g''}$.

Then, in \cite[Theorem 1.(d)]{PZ}, it is asserted that, given a Lie
grading \eqref{eq:Liegr}, the set $G$ embeds in an abelian semigroup
so that the following property holds:
\begin{itemize}
\item[(P)]
For any $g,g',g''\in G$ with $0\ne [\calL_g,\calL_{g'}]\subseteq
\calL_{g''}$, $g+g'=g''$ holds in the semigroup.
\end{itemize}

The purpose of this note is to give a counterexample to this
assertion. The problem in the proof of \cite[Theorem 1.(d)]{PZ} lies
in rule III of \cite[page 104]{PZ}.

\smallskip

Let $V$ be a nine dimensional vector space over a field $k$ with a
fixed basis $\{ a,b_1,b_2,b_3,c_1,c_2,c_3,d_1,d_2\}$, and consider
the endomorphisms $x$, $y$ and $z$ of $V$ whose action on the basic
elements is given in the following diagram
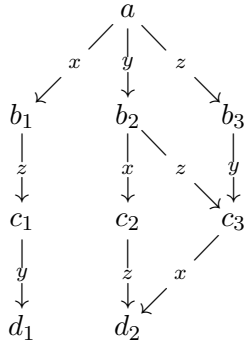
\[
\beginpicture
\setcoordinatesystem units <2pt,2pt>
 \setplotarea x from -40 to 40, y from 0 to 60
 \put{$a$} at 0 60
 \put{$b_1$} at -20 40
 \put{$b_2$} at 0 40
 \put{$b_3$} at 20 40
 \put{$c_1$} at -20 20
 \put{$c_2$} at 0 20
 \put{$c_3$} at 20 20
 \put{$d_1$} at -20 0
 \put{$d_2$} at 0 0
 \multiput{$\diagup$} at -5 55 15 15 /
 \multiput{$\swarrow$} at -15 45 5 5 /
 \multiput{$\diagdown$} at 5 55 5 35 /
 \multiput{$\searrow$} at 15 45 15 25 /
 \multiput{$\vert$} at 0 54 -20 34 0 34 20 34 -20 14 0 14 /
 \multiput{$\downarrow$} at 0 46 -20 26 0 26 20 26 -20 6 0 6 /
 \multiput{${\scriptstyle x}$} at -10 50 0 30 10 10 /
 \multiput{${\scriptstyle y}$} at 0 50 20 30 -20 10 /
 \multiput{${\scriptstyle z}$} at 10 50 -20 30 10 30 0 10 /
\endpicture
\]
(Thus, for instance, $x(a)=b_1$, $x(b_2)=c_2$, $x(c_3)=d_2$ and $x$
annihilates all the other basic elements.)

The associative subalgebra of $\End_k(V)$ generated by these three
endomorphisms is
\[
A=\espan{x,y,z,xy,xz,zx,yz,zy,yzx,xyz}.
\]
Note that $yx=0=x^2=y^2=z^2$, $xyz=xzy=zxy$, $xAx=yAy=zAz=0$, and
$A^4=0$. The elements in the spanning set of $A$ given above are
linearly independent, for if
\begin{multline*}
\alpha_1 x+\alpha_2y+\alpha_3 z+\alpha_4 xy+\alpha_5 xz +\alpha_6 zx
+\alpha_7yz\\
 +\alpha_8 zy+\alpha_9 yzx+\alpha_{10}xyz=0,
\end{multline*}
for some scalars $\alpha_i\in k$, this linear combination applied to
$a$ gives
$\alpha_1=\alpha_2=\alpha_3=\alpha_4=\alpha_6=\alpha_9=\alpha_{10}=0$,
as well as $\alpha_7+\alpha_8=0$. Now applied to $b_2$ gives
$\alpha_5=0$, and to $b_1$ $\alpha_7=0$. Therefore, the dimension of
$A$ is exactly $10$.

Note that
\[
\begin{split}
[[x,y],z]&=xyz-yxz-zxy+zyx=xyz-zxy=0\quad \text{(as $yx=0$),}\\[4pt]
[[y,z],x]&=yzx-zyx-xyz+xzy\\
       &=yzx-0-(xyz-xzy)=yzx\quad \text{(as $xyz=xzy$),}\\[4pt]
[[z,x],y]&=zxy-xzy-yzx+yxz=-yzx.
\end{split}
\]
so the Lie subalgebra of $\End_k(V)$ generated by $x$, $y$ and $z$
is
\[
\frg=\espan{x,y,z,[x,y]=xy, [x,z],[y,z],[[y,z],x]},
\]
which is a seven dimensional nilpotent Lie algebra.

Now, the Lie algebra
\[
\calL=\frg\oplus V\quad \text{(semidirect sum)}
\]
is a nilpotent Lie algebra, and its basis
\[
B=\{x,y,z,[x,y],
[x,z],[y,z],[[y,z],x],a,b_1,b_2,b_3,c_1,c_2,c_3,d_1,d_2\}
\]
satisfies that the bracket of any two elements in $B$ is either $0$
or a scalar multiple of another basic element. Hence, this basis
gives a Lie grading:
\begin{equation}\label{eq:myLiegr}
\calL=\oplus_{u\in B}\calL_u,
\end{equation}
where $\calL_u=ku$ for any $u\in B$.

However, $B$ is not contained in any grading abelian semigroup
satisfying property (P) above, because
\[
[y,[z,[x,a]]]=d_1,\qquad\text{while}\qquad [x,[y,[z,a]]]=d_2,
\]
and $d_1$ and $d_2$ are in different homogeneous components in
\eqref{eq:myLiegr}. If $B$ were contained in an abelian semigroup
satisfying (P), with addition denoted by $\boxplus$, then
\[
d_1=y\boxplus z\boxplus x\boxplus a=x\boxplus y\boxplus z\boxplus
a=d_2
\]
would hold, a contradiction.

\end{document}